\documentclass[11pt,a4paper]{article}
\usepackage[T1]{fontenc}
\usepackage[utf8]{inputenc}
\usepackage{lmodern}
\usepackage{amsmath,amssymb,bm}
\usepackage{booktabs,multirow}
\usepackage{graphicx}
\usepackage[margin=2.5cm]{geometry}
\usepackage{xcolor}
\usepackage {cite}
\usepackage{tcolorbox}
\usepackage{listings}
\usepackage{fancyhdr}
\usepackage{caption,subcaption}
\usepackage[colorlinks=true,linkcolor=blue!60!black,citecolor=blue!60!black,urlcolor=blue!60!black]{hyperref}

\definecolor{titlecol}{RGB}{0,96,96}
\definecolor{boxcol}{RGB}{20,40,110}
\usepackage{sectsty}
\allsectionsfont{\color{titlecol}}
\tcbset{modelbox/.style={colback=white, colframe=boxcol, fonttitle=\bfseries, boxrule=0.6pt, arc=2pt, left=4pt, right=4pt, top=2pt, bottom=2pt}}
\newcommand{\code}[1]{\texttt{#1}}
\newcommand{\BFH}{\textsc{BioFluxHDG}}
\newcommand{\BNF}{\textsc{BioNetFlux}}
\newcommand{\jj}{\bm{j}}
\newcommand{\psivec}{\bm{\psi}}
\newcommand{\nn}{\bm{n}}
\newcommand{\PP}{\mathbb{P}}
\newcommand{\Th}{\mathcal{T}_h}
\newcommand{\Eh}{\mathcal{E}_h}

\begin{document}

\begin{titlepage}
\thispagestyle{empty}
\centering
\vspace*{1.2cm}
\includegraphics[width=0.5\textwidth]{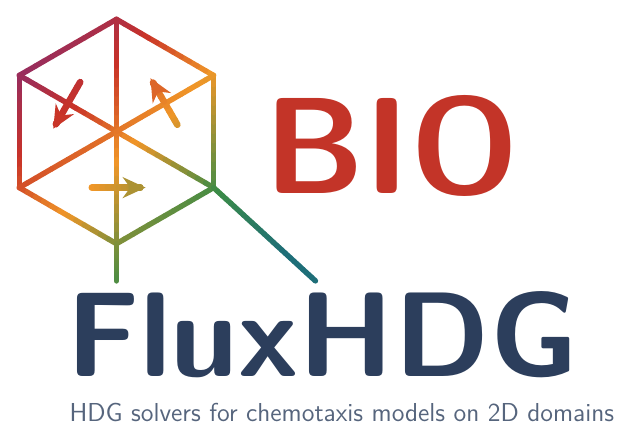}\\[1.6cm]
{\Huge\bfseries BioFluxHDG}\\[0.6cm]
{\Large A MATLAB Framework
	for Reaction--Diffusion--Chemotaxis\\[2pt]
	 Simulations on Two-Dimensional Domains}
	\\[2.2cm]
{\large Silvia Bertoluzza \quad Micol Pennacchio}\\[0.5cm]
{Istituto di Matematica Applicata e Tecnologie Informatiche\\
``E. Magenes'' --- CNR, Pavia, Italy}\\[2.2cm]
{\large November 2026}\\[1.8cm]
{\small\itshape Hybridizable Discontinuous Galerkin methods $\cdot$
	Model-driven assembly\\[3pt]
Keller--Segel chemotaxis $\cdot$ Multiple-sclerosis lesion patterns}
\vfill
\noindent\makebox[\textwidth]{%
  \includegraphics[height=1.35cm]{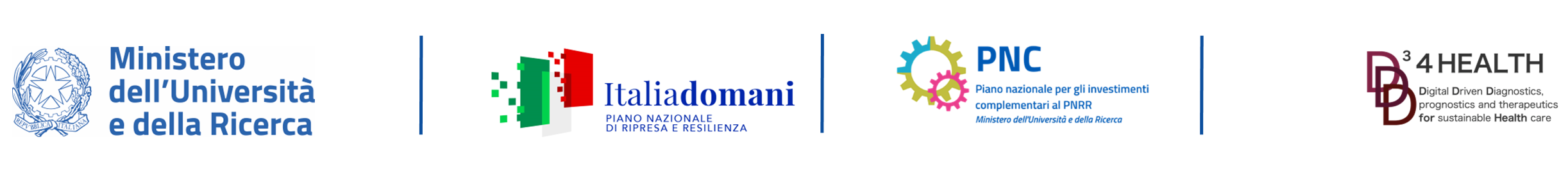}\hspace{0.6cm}%
  \includegraphics[height=0.95cm]{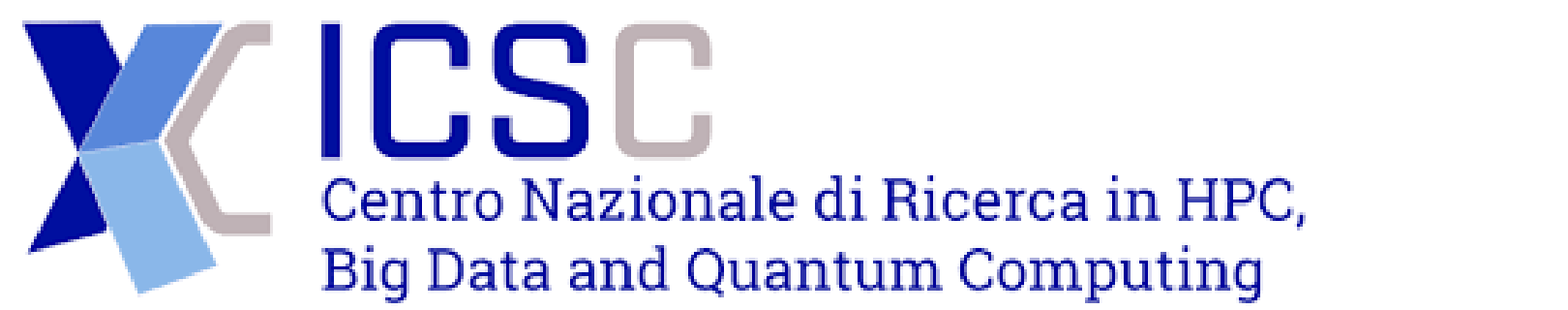}}
\end{titlepage}

\setcounter{page}{1}
\begin{abstract}
We describe \BFH{}, a MATLAB framework for the simulation of
reaction--diffusion--chemotaxis systems on two-dimensional polygonal domains by
hybridizable discontinuous Galerkin (HDG) methods. A model is described as a list
of diffusive species, each with a diffusion coefficient, chemotactic drifts and a
reaction term, and of local species obeying pointwise ODEs; the element
matrices, the linearization and the static condensation are generated from this
description. The code supports the HDG schemes analyzed for the Keller--Segel
system and for a Keller--Segel model of multiple sclerosis, Neumann, Dirichlet
and Robin boundary conditions, implicit Euler and BDF2 time stepping with Picard
or Newton linearization, an upwind stabilization of the chemotactic transport and
a positivity limiter. We present the models, the discretization, the software
architecture and numerical results: convergence tests, the finite-time blow-up of
the Keller--Segel system, and demyelination patterns on a curved domain.
\end{abstract}
\tableofcontents
\clearpage

\section{Introduction}
\BFH{} is a MATLAB framework for the numerical simulation of coupled systems of
reaction--diffusion--chemotaxis equations on two-dimensional polygonal domains.
It is the two-dimensional companion of \BNF{} \cite{bionetflux}, which addresses
the same class of models on one-dimensional networks, and it collects, in a single
model-driven code, the solvers developed for the analysis of the hybridizable
discontinuous Galerkin (HDG) discretization of the Keller--Segel system
\cite{kshdg} and for the simulation of demyelination patterns in a
Keller--Segel-type model of multiple sclerosis \cite{mshdg}.

The framework rests on three ideas.
\begin{enumerate}
\item \textbf{HDG spatial discretization} \cite{cockburn2009,cockburn2016}: every
  species that diffuses is approximated by a scalar, a flux and a trace on the mesh
  skeleton; after static condensation only the traces are globally coupled, and the
  element unknowns are recovered by independent local solves. The pairs
  $\mathrm{HDG}(k',k)$, $k'\in\{k-1,k\}$, $k=1,2$, analyzed in \cite{kshdg} are
  available, including the reduced-degree trace $\widehat u_h\in\PP_{k-1}$ with the
  projected stabilization.
\item \textbf{Model-driven assembly.} A model is a list of \emph{species}: a diffusive
  species carries a diffusion coefficient, optional drifts along the gradient of
  other species (chemotaxis) and a reaction term; a \emph{local} species obeys a
  pointwise ODE and has neither flux nor trace. The element matrices, the
  linearization and the static condensation are generated from this description,
  so that a new model is a short MATLAB file (Section~\ref{sec:architecture}).
\item \textbf{Robust time stepping for aggregation-dominated regimes.} Implicit
  Euler or BDF2 in time; Picard or Newton treatment of the nonlinear terms; an
  upwind stabilization of the chemotactic transport \cite{cockburn_dong2009,nguyen2009b,egger2010},
  optionally combined with a Zhang--Shu positivity limiter \cite{zhangshu2010}, which allows the
  simulation of the finite-time blow-up of the Keller--Segel system while
  preserving the total mass and the positivity of the nodal values.
\end{enumerate}

Key capabilities include: unstructured triangular meshes from the MATLAB PDE
Toolbox geometries; Neumann, Dirichlet and Robin boundary conditions assigned
edge by edge and independently for each species; symbolic generation of
manufactured solutions (sources, fluxes and boundary data) directly from the model
description, for convergence studies; serial or \code{parfor} element loops.
Two models ship with the code: the classical two-equation Keller--Segel system
and the three-equation Keller--Segel model for multiple sclerosis of
\cite{lombardo2017}. The code is available from the authors upon request.

The remainder of this report presents the mathematical models
(Section~\ref{sec:models}), the HDG discretization (Section~\ref{sec:hdg}), the
software architecture (Section~\ref{sec:architecture}), a gallery of results
(Section~\ref{sec:results}) and a getting-started guide
(Section~\ref{sec:start}).

\section{Mathematical Models}\label{sec:models}

\subsection{Species and Equations}
Let $\Omega\subset\mathbb{R}^2$ be a polygonal domain and $T>0$. A model is a
system of $N$ species $u_1,\dots,u_N$. A \emph{diffusive} species satisfies
\begin{equation}\label{eq:diffusive}
\partial_t u_i + \nabla\cdot\jj_i + r_i(u) = f_i, \qquad
\jj_i = -D_i\nabla u_i + \sum_{j=1}^{N} \chi_{ij}\, g_{ij}(u)\,\nabla u_j ,
\end{equation}
where $D_i>0$ is the diffusion coefficient, $\chi_{ij}$ is the chemotactic
sensitivity of species $i$ to the gradient of species $j$, with a possibly
nonlinear factor $g_{ij}(u)$ ($\chi_{ij}=0$ when species $i$ is not attracted by
species $j$; we call \emph{drifts} of species $i$ the pairs $(i,j)$ with
$\chi_{ij}\neq0$), $r_i$ is the reaction term and $f_i$ a source. A \emph{local} species satisfies the pointwise ODE
\begin{equation}\label{eq:local}
\partial_t u_i + r_i(u) = f_i .
\end{equation}
All reaction and drift functions may depend on all species.

\subsection{Keller--Segel Chemotaxis System (2 Equations)}
The classical Keller--Segel model \cite{keller1970} describes a cell density
$u\ge 0$ and a chemoattractant concentration $\varphi\ge0$:
\begin{tcolorbox}[modelbox, title=Keller--Segel system]
\begin{align}
\partial_t u + \nabla\cdot\jj &= f_u, & \jj &= -\mu\nabla u + \chi\, u\,\nabla\varphi, \label{eq:ks1}\\
\partial_t\varphi + \nabla\cdot\psivec + a\varphi - b u &= f_\varphi, & \psivec &= -\nu\nabla\varphi, \label{eq:ks2}
\end{align}
\end{tcolorbox}
\noindent
where $\mu,\nu>0$ are the diffusion coefficients, $\chi>0$ the chemotactic
sensitivity, $b>0$ the production rate of the chemoattractant by the cells and
$a\ge0$ its decay rate. With a constant sensitivity the solution may blow up in
finite time when the initial mass exceeds $8\pi\mu\nu/(b\chi)$ \cite{horstmann2003};
reproducing this regime numerically is the subject of Section~\ref{sec:blowup}.
In the terminology of \eqref{eq:diffusive}: $u$ and $\varphi$ are diffusive with
$D_u=\mu$ and $D_\varphi=\nu$, $u$ drifts towards $\nabla\varphi$ with
$\chi_{u\varphi}=\chi$ and $g(u)=u$, $r_u=0$ and $r_\varphi=a\varphi-bu$.

\subsection{Keller--Segel Model of Multiple Sclerosis (3 Equations)}
The model of \cite{lombardo2017} describes the formation of demyelination
lesions: $u$ is the density of activated macrophages, $\varphi$ the
concentration of cytokines (the chemoattractant) and $d$ the density of
destroyed oligodendrocytes:
\begin{tcolorbox}[modelbox, title=Keller--Segel MS system]
\begin{align}
\partial_t u + \nabla\cdot\jj - u + u^{\alpha} &= f_u, & \jj &= -\mu\nabla u + \widetilde\chi\, g(u)\,\nabla\varphi, \label{eq:ms1}\\
\partial_t\varphi + \nabla\cdot\psivec + a\varphi - b u - \delta d &= f_\varphi, & \psivec &= -\nu\nabla\varphi, \label{eq:ms2}\\
\partial_t d + h(u)\,(d-1) &= f_d, \label{eq:ms3}
\end{align}
with $g(u) = \dfrac{u}{1+u}$, \quad $h(u)=\dfrac{r\,u^2}{1+u}$.
\end{tcolorbox}
\noindent
The chemotactic sensitivity saturates for large densities ($g$), the reaction
$-u+u^\alpha$ ($\alpha=2$ by default) models activation and saturation of the
macrophages, and $d$ is a local species: destroyed oligodendrocytes do not move.
In the terminology of \eqref{eq:diffusive}--\eqref{eq:local}: $D_u=\mu$,
$D_\varphi=\nu$, $\chi_{u\varphi}=\widetilde\chi$ with the saturating factor $g(u)$,
$r_u=-u+u^\alpha$, $r_\varphi=a\varphi-bu-\delta d$ and $r_d=h(u)(d-1)$.
Depending on $(\nu,\widetilde\chi)$ the homogeneous equilibrium
$(u,\varphi,d)=(1,(b+\delta)/a,1)$ is Turing-unstable and the system develops
hexagonal, striped or spotted patterns; the three presets
$(\nu,\tilde\chi)=(0.8,7.32)$, $(0.8,8.33)$, $(0.02,2.62)$ of the code correspond
to these regimes, and are taken from \cite[Fig.~3]{bisi2025}. Note that in \cite{mshdg} the diffusion coefficient of
$u$ is set to one and that of $\varphi$ is denoted by $\mu$, so that the presets
read $(\mu,\chi)=(0.8,7.32)$, $(0.8,8.33)$ and $(0.02,2.62)$ there.

\subsection{Boundary Conditions}
For every diffusive species, independently, each boundary edge carries one of
\begin{align*}
\text{Neumann:}\quad & -\jj_i\cdot\nn = g^N_i, &
\text{Dirichlet:}\quad & u_i = g^D_i, &
\text{Robin:}\quad & -\jj_i\cdot\nn + u_i/\varepsilon_i = g^R_i ,
\end{align*}
with $\nn$ the outward unit normal and $\varepsilon_i>0$. The default is the
homogeneous Neumann (zero flux) condition on the whole boundary. Edges are
tagged by a user function of the edge midpoint (one tag per side of a rectangle,
or by distance from given points) or by listing points that the edge must
contain; the latter is meant for the small inlets and membranes of
organ-on-chip geometries. Local species have no boundary condition.

\section{HDG Spatial Discretization}\label{sec:hdg}
Let $\Th$ be a shape-regular triangulation of $\Omega$ and $\Eh$ its skeleton.
For $m\ge0$, $\PP_m(\Th)$ and $\PP_m(\Eh)$ denote the discontinuous piecewise
polynomials of degree $m$ on the elements and on the edges. The scheme
$\mathrm{HDG}(k',k)$, $k'\in\{k-1,k\}$, $k\in\{1,2\}$, approximates a diffusive
species $i$ by
\[
u_{i,h}\in\PP_k(\Th),\qquad \jj_{i,h}\in\PP_{k_i}(\Th)^2,\qquad \widehat u_{i,h}\in\PP_{k_i}(\Eh),
\]
where $k_i\in\{k',k\}$ is chosen species by species (the \emph{flux degree} of
the model description). In the Keller--Segel model of Section~2.2 the
	cellular flux may use either $k'$, as in \cite{chen2019,kshdg}, or $k$. The
	equal-order choice $k_i=k$, which gives order $k+1$ also for the cellular flux,
	is the one analyzed in \cite{mshdg} for the multiple-sclerosis model. Local
species are approximated in $\PP_k(\Th)$ with no flux and no trace.

\subsection{Local Problem and Numerical Fluxes}
On every element $K$, with $\alpha_i=1/D_i$, the mixed form of
\eqref{eq:diffusive} is tested against $(w,\bm q)\in\PP_k(K)\times\PP_{k_i}(K)^2$:
\begin{align}
(\partial_t u_{i,h}, w)_K - (\jj_{i,h},\nabla w)_K + \langle \widehat\jj_{i,h}\cdot\nn, w\rangle_{\partial K}
 + (r_i(u_h), w)_K &= (f_i, w)_K, \label{eq:loc1}\\
(\alpha_i\jj_{i,h},\bm q)_K - (u_{i,h},\nabla\cdot\bm q)_K + \langle\widehat u_{i,h},\bm q\cdot\nn\rangle_{\partial K}
 + \sum_{j=1}^{N} \frac{\chi_{ij}}{D_iD_j}\big(g_{ij}(u_h)\,\psivec_{j,h},\bm q\big)_K &= 0, \label{eq:loc2}
\end{align}
where the drift has been expressed through the flux $\psivec_{j,h}=-D_j\nabla u_{j,h}$
of the attracting species. The numerical normal flux is
\begin{equation}\label{eq:numflux}
\widehat\jj_{i,h}\cdot\nn_K = \jj_{i,h}\cdot\nn_K + \gamma_i\, h^{-\kappa_i}\,
\Pi^\partial_{k_i}\big(u_{i,h}-\widehat u_{i,h}\big), \qquad \kappa_i = k-k_i,
\end{equation}
with $\gamma_i>0$ a stabilization constant and $\Pi^\partial_{k_i}$ the $L^2$
projection onto $\PP_{k_i}(e)$ on each edge. For $k_i=k$ the projection is the
identity. For $k_i=k-1$ two variants are implemented: the trace in
$\PP_{k-1}(\Eh)$ with the projected jump, which is the scheme analyzed in
\cite{kshdg} (option \code{trace\_degree = 'kprime'}), and the trace in
$\PP_k(\Eh)$ without projection, as in the original codes (default).

\paragraph{Transmission conditions.} The traces are single valued on $\Eh$ and
the normal fluxes balance across every interior edge,
$\langle\widehat\jj_{i,h}\cdot\nn,\widehat w\rangle_{\Eh}=0$ for all
$\widehat w\in\PP_{k_i}(\Eh)$, while on boundary edges the same equation carries
the Neumann or Robin datum; Dirichlet values are imposed on the trace.

\subsection{Static Condensation}
Collecting on $K$ the diffusive unknowns $U_K=(u_{i,h},\jj_{i,h})_i$ and the
traces $\widehat U_K$, one time step of \eqref{eq:loc1}--\eqref{eq:loc2} reads
\begin{equation}\label{eq:local_system}
K_2\, U_K = g + B\,\widehat U_K, \qquad
B^{\mathsf T}U_K - \widehat A\,\widehat U_K = g_\partial ,
\end{equation}
where $K_2$ collects the static operators (divergence, normal trace, vector mass),
the mass terms $M/\Delta t$, the linearized reactions and drifts and the
stabilization \eqref{eq:numflux}; $\widehat A$ is the stabilization/Robin mass on
$\partial K$ and $g_\partial$ holds the boundary data. Eliminating $U_K$,
\[
\big(\widehat A - B^{\mathsf T}K_2^{-1}B\big)\,\widehat U_K = B^{\mathsf T}K_2^{-1}g - g_\partial,
\qquad U_K = K_2^{-1}g + K_2^{-1}B\,\widehat U_K .
\]
The element contributions are assembled into a sparse global system in the
traces only, solved by a direct solver, and the element unknowns are recovered
locally. The local blocks are dense and small ($\le 36$ unknowns for
$\mathrm{HDG}(2,2)$ with two species) and $K_2$ is factorized once for both
right-hand sides. The element loop is independent across elements and runs in
\code{parfor} when a pool is open.

\subsection{Linearization of Reactions and Drifts}\label{sec:linearization}
Time derivatives are discretized by the implicit Euler method or by BDF2
(the latter written as an implicit Euler step with $\Delta t_{\mathrm{eff}} = 2\Delta t/3$
and old value $(4u^n-u^{n-1})/3$; the first step is an Euler step). At each step
the nonlinear system is solved by a fixed-point iteration in which, around the
current iterate $u^{(m)}$,
\begin{itemize}
\item a reaction $r_i(u)$ is linearized by Newton with respect to the species for
  which the model provides the derivative $\partial r_i/\partial u_j$, and frozen
  (Picard) with respect to the others;
\item the drift term $g_{ij}(u)\psivec_j$ is linearized either with the flux frozen
  (\emph{Picard}, default): $g_{ij}(u)\psivec_j\approx[g_{ij}(u^{(m)})+g'_{ij}(u^{(m)})(u_i-u_i^{(m)})]\psivec_j^{(m)}$,
  or fully (\emph{Newton}): $g_{ij}(u^{(m)})\psivec_j + g'_{ij}(u^{(m)})\psivec_j^{(m)}(u_i-u_i^{(m)})$,
  which adds a block coupling $\jj_i$ and $\psivec_j$ inside the element;
\item a local species is updated first, element by element, by a Newton step on
  \eqref{eq:local} with the diffusive species at the current iterate; its new
  value enters the reactions of the diffusive species.
\end{itemize}
The iteration stops when the $L^2$ increments of a user-selected set of species
fall below their tolerances (by default all diffusive species, tolerance
$10^{-12}$). For the Keller--Segel system Picard converges in 3 iterations per
step in the convergence tests and in 5--10 in the blow-up regime; Newton in 2
and 4 respectively. Picard is the default because adding an equation to the
model does not require any new Jacobian: the same freezing procedure handles the
new nonlinear terms (see also \cite[Remark~5.1]{kshdg}).

\subsection{Upwind Stabilization and Positivity Limiter}\label{sec:upwind}
In aggregation-dominated regimes the chemotactic transport generates spurious
oscillations and negative densities. As in \cite{kshdg} and \cite{cockburn_dong2009}, where the total flux is likewise taken as an
independent unknown, the stabilization constant of a drifting species in the
numerical flux \eqref{eq:numflux} is split as
\begin{equation}\label{eq:upwind_split}
\gamma_i = \gamma_D + \gamma_C,
\end{equation}
where $\gamma_D>0$ is the diffusive contribution and $\gamma_C$ is the upwind
contribution associated with the chemotactic transport, acting on the inflow
part of $\partial K$; see also \cite{nguyen2009a,nguyen2009b,egger2010} for the
case in which the convective flux is treated separately. For each pair $(i,j)$
with $\chi_{ij}\neq0$, the chemotactic velocity is
$\bm b_{ij} = \chi_{ij}\, v_{ij}(u)\,\nabla u_j = -(\chi_{ij}/D_j)\, v_{ij}(u)\,\psivec_j$,
with $v_{ij}=g_{ij}(u)/u_i$ the velocity factor ($v=1$ for the Keller--Segel
drift $\chi u\nabla\varphi$, $v=1/(1+u)$ for the MS model), so that the inflow
condition $\bm b_{ij}\cdot\nn_K<0$ is equivalent to $\psivec_j\cdot\nn_K>0$.
Accordingly, the convective stabilization parameter is evaluated edgewise as
\begin{equation}\label{eq:upwind}
\gamma_C\big|_e = \sum_{j=1}^{N}
\max\Big(\frac{\chi_{ij}}{D_j}\,\frac{1}{|e|}\int_e v_{ij}(u_h^{(m)})\,
\psivec_{j,h}^{(m)}\cdot\nn_K,\ 0\Big), \qquad e\subset\partial K,
\end{equation}
which for the Keller--Segel system reduces to \cite[eq.~(5.3)]{kshdg},
$\gamma_C|_e=\max\big(\beta\mu\widetilde\chi\,\frac{1}{|e|}\int_e\psivec_h\cdot\nn_K,0\big)$
with $\beta=1/\nu$, $\widetilde\chi=\chi/\mu$. The edge integral is computed by
Simpson's rule on the edge nodal values (option \code{upwind\_mean = 'integral'});
the arithmetic mean of the edge nodal values is available as \code{'nodal'}. The
quantities in \eqref{eq:upwind} are frozen at the current iterate, also with the
Newton linearization.

Optionally, after every time step the nodal values of a chosen species are
corrected on each element by the scaling limiter of Zhang and Shu
\cite{zhangshu2010}: with $\bar u_K$ the cell mean and $m_K$ the nodal minimum,
if $m_K<\epsilon_K=\epsilon\,\bar u_K$ the values are replaced by
$\bar u_K+\theta_K(u_h-\bar u_K)$, $\theta_K=(\bar u_K-\epsilon_K)/(\bar u_K-m_K)$,
so that the nodal minimum becomes $\epsilon_K\ge0$ and $\int_K u_h$ is unchanged.
The traces and the other species are not modified.

\subsection{Manufactured Solutions}
Given exact expressions $u_i(x,y,t)$ for all species, the sources
$f_i=\partial_t u_i+\nabla\cdot\jj_i+r_i(u)$ (or $\partial_t u_i+r_i(u)$ for local
species), the exact fluxes and the Dirichlet/Robin data are generated
symbolically from the model description (Symbolic Math Toolbox), since the
reaction and drift functions are ordinary MATLAB handles evaluated on symbolic
variables. This is how the convergence tables of Section~\ref{sec:results} are
produced and how a user-defined model can be verified in a few lines.

\section{Software Architecture}\label{sec:architecture}
\BFH{} is a plain MATLAB code (R2025a or later; Parallel Computing Toolbox
optional, Symbolic Math Toolbox for manufactured solutions, PDE Toolbox for the
mesh generator). It is organized in two packages and a models folder:
\begin{tcolorbox}[modelbox, title=Package layout]
\begin{verbatim}
BioFluxHDG/
+-- +hdg/      generic HDG layer, no physics: mesh and trace numbering,
|              Lagrange bases and quadrature, boundary tags per species,
|              Zhang-Shu limiter, field plots
+-- +chemo/    model-driven solver: Model (species, drifts, reactions),
|              config, Solver (layout, traces, initial data, time loop,
|              nonlinear iterations, condensation), element kernel,
|              static matrices, sources, L2 errors, manufactured solutions
+-- models/    ks.m (Keller-Segel), ks_ms.m (multiple sclerosis)
+-- examples/  ks_convergence, ks_blowup, ms_convergence, ms_pattern,
|              ms_pattern_4leaf_run
+-- tests/     regression against the original codes
+-- meshes/    PDE Toolbox geometries
\end{verbatim}
\end{tcolorbox}

\paragraph{Model description.} A model file builds a \code{chemo.Model}:
\begin{lstlisting}[caption={The Keller--Segel model (\code{models/ks.m}, abridged).}]
model = chemo.Model('Keller-Segel', params);
model = model.add_species('u',   'diffusive', 'D', mu, 'flux_degree', 'kprime', 'stab', Ks);
model = model.add_species('phi', 'diffusive', 'D', nu, 'stab', 1, ...
            'reaction',  @(S) a*S.phi - b*S.u, ...
            'dreaction', struct('phi', @(S) a + 0*S.phi, 'u', @(S) -b + 0*S.u));
model = model.add_drift('u', 'phi', chi);          % g(S) = S.u, velocity factor 1
\end{lstlisting}
Species values are passed to the reaction and drift functions as a struct
\code{S} with one field per species (vectors at the quadrature points). A local
species is declared with kind \code{'local'}; a nonlinear drift factor with
\code{'f'}, its derivative with \code{'df'} and the velocity factor of
\eqref{eq:upwind} with \code{'upwind\_factor'}. Terms with a derivative are
treated by Newton, the others frozen, as described in
Section~\ref{sec:linearization}. The flux degree (\code{'k'} or
\code{'kprime'}) and the stabilization constant are properties of the species;
the element pair, the trace degree and all numerical options are properties of
the configuration (\code{chemo.config}).

\paragraph{Data flow.} \code{chemo.Solver(model, config, problem)} receives the
model, the configuration and the \emph{problem} data (initial data, sources,
exact fluxes and boundary data per species, either written by the user or
generated by \code{chemo.manufactured}). \code{setup(h)} builds the mesh, the
layout of the element unknowns and of the traces (one global numbering per
diffusive species, with the trace degree of that species), the static element
matrices and the boundary tags. \code{run\_time\_integration()} advances in
time: at every step, sources and boundary data are assembled at $t^{n+1}$ and
the nonlinear iterations call, for every element in \code{parfor}, the kernel
\code{chemo.condense\_element}, which updates the local species, adds the
linearized reactions, drifts and the (upwind) stabilization to the static
operator and returns the condensed blocks; the global trace system is assembled
from triplets and solved; the element unknowns are recovered. The solver
returns the traces, the element unknowns of the diffusive and local species at
the stored time levels, the iteration counts and the limiter statistics.

\paragraph{Configuration.} \code{chemo.config('Name', value, ...)} collects the
numerical options: element pair (\code{'P0P1'}, \code{'P1P2'}, \code{'P1P1'},
\code{'P2P2'}), trace degree, mesh sizes, final time and time step, time scheme
(\code{'BE'}, \code{'BDF2'}), linearization (\code{'picard'}, \code{'newton'}),
tolerance and stopping species, upwind and edge mean, limiter, random
perturbation of the initial data (nodal or smooth Fourier field with a fixed
seed, for pattern formation), positivity check, storage frequency, parallel
execution, domain and boundary data. Physical parameters belong to the model
file.

\paragraph{Validation.} The framework was built from two validated codes: the
Keller--Segel solver of \cite{kshdg} and the multiple-sclerosis solver of
\cite{mshdg}. Twenty regression cases (all element pairs; Neumann, Dirichlet and
Robin edges; upwind on and off; reduced-degree trace; non-unit physical
parameters; Newton and BDF2; blow-up data with and without limiter; pattern runs
with perturbed initial data) reproduce the reference results of the original
codes to $10^{-14}$ in relative norm with identical iteration counts; parallel
and serial runs coincide bit for bit.

\section{Simulation Results}\label{sec:results}

\subsection{Convergence for the Keller--Segel System}
We reproduce Example~1 of \cite{kshdg}: $\Omega=[0,1]^2$, exact solution
$u=\cos t\,\sin\pi x\cos\pi y$, $\varphi=\sin t\,\cos\pi x\sin\pi y$, all
parameters equal to~1, non-homogeneous Neumann data on the whole boundary,
unstructured meshes of maximum edge length $h=2^{-1},\dots,2^{-4}$, implicit
Euler with $\Delta t=2\cdot10^{-4}$, $T=1$, $\gamma_u=\gamma_\varphi=1$, no
upwind, Picard iterations (3 per step). For $\mathrm{HDG}(0,1)$ and
$\mathrm{HDG}(1,2)$ the trace of $u$ is in $\PP_{k-1}$ (\code{trace\_degree =
'kprime'}). Table~\ref{tab:ks_conv} reports the $L^2(\Omega)$ errors at $t=T$ and
the computed orders: the scalar variables and the chemical flux converge with
order $k+1$ and the cellular flux with order $k'+1$, in agreement with the
analysis of \cite{kshdg}.

\begin{table}[htbp]
\centering
\caption{Keller--Segel, Example 1 of \cite{kshdg}: $L^2(\Omega)$ errors at $T=1$
and computed convergence orders (\code{examples/ks\_convergence.m}).}
\label{tab:ks_conv}
\footnotesize
\setlength{\tabcolsep}{3.5pt}
\begin{tabular}{ll lc lc lc lc}
\toprule
 & & \multicolumn{2}{c}{$\|u-u_h\|$} & \multicolumn{2}{c}{$\|\varphi-\varphi_h\|$}
   & \multicolumn{2}{c}{$\|\jj-\jj_h\|$} & \multicolumn{2}{c}{$\|\psivec-\psivec_h\|$} \\
\cmidrule(lr){3-4}\cmidrule(lr){5-6}\cmidrule(lr){7-8}\cmidrule(lr){9-10}
 & $h$ & Error & Order & Error & Order & Error & Order & Error & Order \\
\midrule
\multirow{4}{*}{HDG(0,1)}
 & $2^{-1}$ & 2.47e-1 & --   & 8.47e-2 & --   & 6.70e-1 & --   & 1.65e-1 & --   \\
 & $2^{-2}$ & 6.67e-2 & 1.89 & 2.71e-2 & 1.64 & 3.70e-1 & 0.86 & 5.17e-2 & 1.67 \\
 & $2^{-3}$ & 1.59e-2 & 2.07 & 6.89e-3 & 1.98 & 1.75e-1 & 1.08 & 1.30e-2 & 2.00 \\
 & $2^{-4}$ & 3.89e-3 & 2.03 & 1.72e-3 & 2.00 & 8.50e-2 & 1.04 & 3.18e-3 & 2.03 \\
\midrule
\multirow{4}{*}{HDG(1,1)}
 & $2^{-1}$ & 7.36e-2 & --   & 8.17e-2 & --   & 1.70e-1 & --   & 1.49e-1 & --   \\
 & $2^{-2}$ & 2.36e-2 & 1.64 & 2.65e-2 & 1.62 & 5.34e-2 & 1.67 & 4.82e-2 & 1.63 \\
 & $2^{-3}$ & 5.66e-3 & 2.06 & 6.75e-3 & 1.97 & 1.37e-2 & 1.97 & 1.22e-2 & 1.99 \\
 & $2^{-4}$ & 1.41e-3 & 2.01 & 1.69e-3 & 2.00 & 3.44e-3 & 1.99 & 2.98e-3 & 2.03 \\
\midrule
\multirow{4}{*}{HDG(1,2)}
 & $2^{-1}$ & 3.67e-2 & --   & 1.34e-2 & --   & 1.58e-1 & --   & 2.89e-2 & --   \\
 & $2^{-2}$ & 5.97e-3 & 2.62 & 2.27e-3 & 2.56 & 4.29e-2 & 1.88 & 4.15e-3 & 2.80 \\
 & $2^{-3}$ & 7.14e-4 & 3.06 & 2.65e-4 & 3.10 & 1.06e-2 & 2.02 & 4.81e-4 & 3.11 \\
 & $2^{-4}$ & 8.84e-5 & 3.01 & 3.18e-5 & 3.06 & 2.59e-3 & 2.03 & 5.75e-5 & 3.07 \\
\midrule
\multirow{4}{*}{HDG(2,2)}
 & $2^{-1}$ & 1.28e-2 & --   & 1.34e-2 & --   & 3.09e-2 & --   & 2.86e-2 & --   \\
 & $2^{-2}$ & 2.57e-3 & 2.31 & 2.27e-3 & 2.56 & 5.87e-3 & 2.40 & 4.13e-3 & 2.79 \\
 & $2^{-3}$ & 3.15e-4 & 3.03 & 2.65e-4 & 3.10 & 7.25e-4 & 3.02 & 4.81e-4 & 3.10 \\
 & $2^{-4}$ & 3.66e-5 & 3.11 & 3.18e-5 & 3.06 & 8.57e-5 & 3.08 & 5.75e-5 & 3.06 \\
\bottomrule
\end{tabular}
\end{table}

\subsection{Blow-up of the Keller--Segel System}\label{sec:blowup}
Example~2 of \cite{kshdg}: $\Omega=[-\tfrac12,\tfrac12]^2$, homogeneous Neumann
conditions, Gaussian initial data $u_0=840\,e^{-84 r^2}$, $\varphi_0=420\,e^{-42 r^2}$
\cite{epshteyn2009}, whose mass $10\pi$ exceeds the critical value $8\pi$,
$\mathrm{HDG}(2,2)$ on a mesh with $h=1/70$, $\Delta t=10^{-5}$,
$T=1.2\cdot10^{-4}$ (\code{examples/ks\_blowup.m}). Without upwind the computed
density develops spurious oscillations and negative nodal values at
$t=6\cdot10^{-5}$, where the run stops by the positivity check. With the
upwind stabilization \eqref{eq:upwind} ($\gamma_D=1$, integral edge mean) the
density concentrates at the origin, the nodal values of $u_h$ stay positive at
all time steps and the total mass is preserved to round-off,
$\max_n|m(t_n)-m(0)|/m(0)\approx 6\cdot10^{-16}$. 
The same holds when the time
step $\Delta t$ is halved. The code also allows the use of the Zhang--Shu
limiter of Section~\ref{sec:upwind}, through the option
\code{limiter = \{'u'\}}, which was not used in \cite{kshdg}. The density
profiles of this test are shown in \cite{kshdg}.

\subsection{Demyelination Patterns in the Multiple-Sclerosis Model}
\label{sec:patterns}
\begin{figure}[htbp]
\centering
\includegraphics[width=\textwidth]{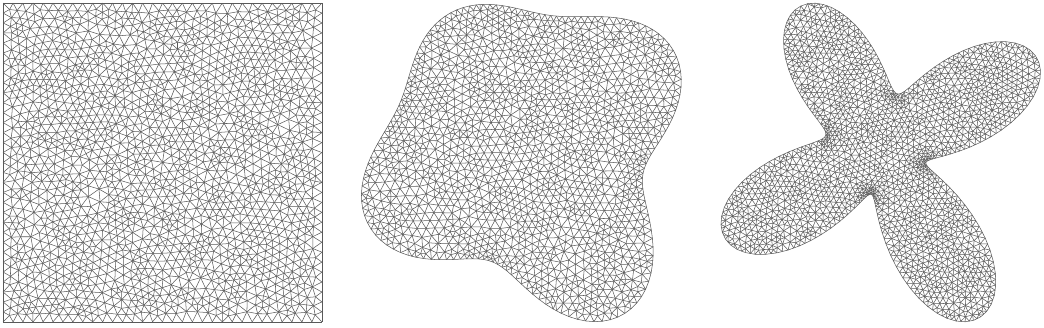}
\caption{The three geometries available for the pattern runs and their meshes
with $h\approx0.59$: square $[0,6\pi]^2$ (\code{gpi6}, 3402 triangles), mildly
four-lobed domain (\code{g4leaf}, 2963) and four-leaf clover (\code{g4leaf055},
3556).}
\label{fig:mesh4leaf}
\end{figure}

\begin{figure}[htbp]
\centering
\includegraphics[width=\textwidth]{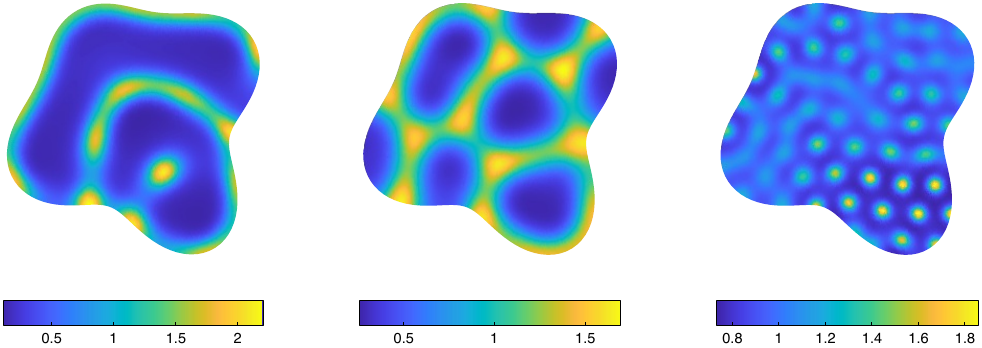}
\caption{MS model on the four-lobed domain: $u_h$ for the presets stripes
($t=500$), hexagons ($t=1000$) and isolated spots ($t=1000$), computed
by \code{examples/ms\_pattern\_4leaf\_run.m}.}
\label{fig:patterns}
\end{figure}

The code ships three geometries of size comparable to $6\pi$ for the pattern runs
(Figure~\ref{fig:mesh4leaf}): the square $[0,6\pi]^2$ of \cite{lombardo2017},
a domain with a smooth, mildly four-lobed boundary (\code{meshes/g4leaf.mat})
and a four-leaf clover (\code{meshes/g4leaf055.mat}), all meshed by Delaunay
triangulation with $h\approx0.59$. The patterns on the square and on the four-leaf
clover are reported in \cite{mshdg}; here we show the mildly four-lobed domain
(2963 triangles). With
homogeneous Neumann conditions, the homogeneous equilibrium
$(u,\varphi,d)=(1,(b+\delta)/a,1)$ is perturbed at $t=0$ by a random nodal
field of amplitude $10^{-2}$ (fixed seed, so that runs are reproducible) and
the system is integrated with $\mathrm{HDG}(2,2)$, Newton linearization and
BDF2 with $\Delta t=0.1$ (\code{examples/ms\_pattern\_4leaf\_run.m}; about
$3$ Newton iterations per step). Figure~\ref{fig:patterns} shows the macrophage
density $u_h$ for the three presets, stripes, hexagons and isolated
spots, selected by $(\nu,\tilde\chi)$ according to the linear stability
analysis of \cite{lombardo2017}  and to the pattern analysis of \cite{bisi2025}. The stripes appear at $t\approx100$ and then
adapt to the geometry: they meet the boundary at a right angle, as imposed
by the no-flux condition, or run along it; the hexagonal lattice is
stationary from $t\approx400$ on; the spots preset ($\nu=0.02$) has the
slowest growth and the lattice of spots, visible at $t=1000$, is still growing
in amplitude.

\clearpage
\section{Getting Started}\label{sec:start}
The listings below show how to install the code, how to run a convergence
test with a manufactured solution, how to reproduce the blow-up test of
Section~\ref{sec:blowup}, and how to define a new model.
\begin{lstlisting}[caption={Installation and path.}]
% unzip / clone BioFluxHDG, then in MATLAB:
cd BioFluxHDG
setpath                       % adds the package root, models/ and meshes/
\end{lstlisting}

\begin{lstlisting}[caption={Keller--Segel convergence test with a manufactured solution.}]
[model, data] = ks('mu', 1, 'nu', 1, 'chi', 1, 'Ks', 1);
cfg = chemo.config('fem', 'P1P2', 'trace_degree', 'kprime', 'h_values', [1/4 1/8 1/16], ...
                   'T_final', 1, 'dt', 2e-4, 'upwind', false, 'stop_species', {'phi', 1});
[problem, cfg] = chemo.manufactured(model, data.exact, cfg);   % sources, fluxes, BC data
for h = cfg.mesh.h_values
    s = chemo.Solver(model, cfg, problem);
    s = s.setup(h);
    [hU, Uk, Wk, times, iters] = s.run_time_integration();
    err = chemo.l2_errors(s, Uk, Wk, times);     % err.u, err.phi, err.flux.u, err.flux.phi
    fprintf('h = %.4f  |u-u_h| = %.3e  |phi-phi_h| = %.3e\n', h, err.u(end), err.phi(end));
end
\end{lstlisting}

\begin{lstlisting}[caption={Blow-up with upwind stabilization and limiter.}]
[model, data] = ks('Ks', 1);
dom = struct('bounds', [-0.5 0.5; -0.5 0.5], 'mesh_file', 'g2', 'mesh_var', 'g');
cfg = chemo.config('fem', 'P2P2', 'h_values', 1/70, 'T_final', 1.2e-4, 'dt', 1e-5, ...
                   'upwind', true, 'upwind_mean', 'integral', 'limiter', {'u'}, ...
                   'exact_solution', false, 'domain', dom, 'parallel', true);
s = chemo.Solver(model, cfg, data.epshteyn);  s = s.setup(1/70);
[hU, Uk, Wk, times, iters, info] = s.run_time_integration();
hdg.plot_field(s.field(Uk(:,:,end), [], 'u'), s.mesh.region, 'u_h', times(end), 2);
\end{lstlisting}

\begin{lstlisting}[caption={A user-defined model: Keller--Segel with a local species.}]
function [model, data] = my_model(p)
    model = chemo.Model('KS + damage', p);
    model = model.add_species('u',   'diffusive', 'D', p.mu, ...
                'flux_degree', 'kprime', 'stab', p.Ks);
    model = model.add_species('phi', 'diffusive', 'D', p.nu, 'stab', 1, ...
                'reaction', @(S) p.a*S.phi - p.b*S.u - p.delta*S.w, ...
                'dreaction', struct('phi', @(S) p.a + 0*S.phi, 'u', @(S) -p.b + 0*S.u));
    model = model.add_species('w', 'local', ...
                'reaction', @(S) p.r*S.u.*(S.w - 1), 'dreaction', struct('w', @(S) p.r*S.u));
    model = model.add_drift('u', 'phi', p.chi, 'f', @(S) S.u./(1+S.u), ...
                'df', @(S) 1./(1+S.u).^2, 'upwind_factor', @(S) 1./(1+S.u));
    data.exact = struct('u', @(x,y,t) 1 + 0.5*cos(t).*sin(pi*x).*cos(pi*y), ...
                        'phi', @(x,y,t) sin(t).*cos(pi*x).*sin(pi*y), ...
                        'w', @(x,y,t) (1 + cos(pi*x).*cos(pi*y)).*(t+1));
end
\end{lstlisting}
The last listing is enough to run a convergence test of the new model with
\code{chemo.manufactured}, since sources and boundary data are generated from
the model description.

\section{Conclusions and Outlook}
\BFH{} provides a compact, extensible platform for the simulation of
reaction--diffusion--chemotaxis systems on two-dimensional domains with
hybridizable discontinuous Galerkin methods. The model-driven assembly makes
the addition of species and nonlinear couplings a matter of a few lines, while
the HDG structure keeps the global system small (traces only)
while the element computations are independent and can be carried out in parallel.
The upwind stabilization and the positivity
limiter extend the applicability to aggregation-dominated regimes.

\

Current development directions include:
\begin{itemize}
\item cancer-on-chip models with several cell populations and signalling
  molecules, as in the one-dimensional setting of \BNF{};
\item membrane (Kedem--Katchalsky) conditions between subdomains and the
  coupling with one-dimensional network models;
\item adaptive time stepping driven by the nonlinear iteration count, as in
  \BNF{}, and higher-order time schemes.
\end{itemize}

\section*{Acknowledgements}
\addcontentsline{toc}{section}{Acknowledgements}
In the preparation of this code and report, Claude (Anthropic) was used to
assist in restructuring the original MATLAB implementations, written by the
authors, into the model-driven framework, and in drafting the documentation.
All code and text were reviewed, verified and validated by the authors.

This work was supported by the Italian Ministry of Research, under the
complementary actions to the NRRP ``D34Health -- Digital Driven Diagnostics,
prognostics and therapeutics for sustainable Health care'' Grant \#PNC0000001,
CUP B53C22006100001. 
MP has been partially supported by ICSC -- Centro Nazionale di Ricerca in High
Performance Computing, Big Data, and Quantum Computing funded by the European
Union -- NextGenerationEU plan. The authors are members of the Istituto Nazionale
di Alta Matematica -- Gruppo Nazionale di Calcolo Scientifico.

\end{document}